\documentclass[11pt]{amsart}
\usepackage{amsmath}
\usepackage{amsfonts}
\usepackage{amssymb}
\usepackage{amscd}
\input{epsf.sty}
\catcode `\@=11
\def\numberbysection{\@addtoreset{equation}{section}
         \renewcommand{\theequation}{\thesection.\arabic{equation}}}
\numberbysection
\def\subsubsection{\@startsection{subsubsection}{3}%
  \normalparindent{.5\linespacing\@plus.7\linespacing}{-.5em}%
  {\normalfont\bfseries}}
\newtheorem{thm}{Theorem}[section]
\newtheorem{lem}[thm]{Lemma}
\newtheorem{prop}[thm]{Proposition}
\newtheorem{cor}[thm]{Corollary}

\theoremstyle{definition}

\newtheorem{df}[thm]{Definition}

\newtheorem{rmk}[thm]{Remark}

\def\nn{\nonumber}
\def\vardel{\delta}

\def\In{In}

\def\Sn{\mathbb{S}_n}

\def\Zz{\mathbb{Z}/2\mathbb{Z}}

\def\del{\partial}
\def\s{\sigma}
\def\D{\Delta}

\def\t{\tau}
\def\G{\Gamma}
\def\a{\alpha}

\def\cact{\mathcal{C}act}
\def\Cact{\mathcal{C}act}
\def\Cacti{\mathcal{C}acti}

\def\CWcact{K}
\def\CCCacti{CC_*(Cacti^1)}
\def\CWtwo{K'}

\def\whitevert{v_{white}}

\def\num{{num}}

\def\Zdec{dec^{\pm}}
\def\sdec{dec'}
\def\Fnum{F_{\prec'_v}}
\def\val{val}

\def\color{{clr}}
\def\mk{mk}

\def\lab{{Lab}}

\def\wlbptree{\mathcal{T}_{bp}^{pp,nt}}
\def\swlbptree{\mathcal{T}_{bp}^{pp,nt,\sdec}}

\def\cacti{\mathcal{C}acti}
\def\cact{\mathcal{C}act}

\def\Arc{\mathcal{A}rc}

\begin{document}

\title{A proof of a cyclic version of Deligne's conjecture via Cacti}

\author
[Ralph M.\ Kaufmann]{Ralph M.\ Kaufmann}
\email{kaufmann@math.uconn.edu}

\address{University of Connecticut, Department of Mathematics,
Storrs, CT 06269}

\begin{abstract}
In this note, we show that the normalized Hochschild co--chains of
an associative algebra with a non--degenerate, symmetric,
invariant inner product are an algebra over a chain model of the
framed little discs operad which is given by cacti. In particular,
in this sense they are a BV algebra up to homotopy and the
Hochschild cohomology of such an algebra is a BV algebra whose
induced bracket coincides with Gerstenhaber's bracket. To show
this, we use a cellular chain model for the framed little disc
operad in terms of normalized cacti. This model is given by
tensoring our chain model for the little discs operad in terms of
spineless cacti with natural chain models for $(S^1)^{\times n}$
adapted to cacti.
\end{abstract}

\maketitle

\section*{Introduction}

In this note, we expand our  chain model of the little discs
operad which we gave in terms of spineless cacti to a chain model
for the framed little discs operad in terms of normalized cacti.
Extending the philosophy of \cite{del}, we then show that the
chain model for the framed little discs operad naturally acts on
the normalized Hochschild cochains of a unital associative algebra
with a non--degenerate, symmetric, invariant bi--linear pairing.
In fact, as in \cite{del}, this operation  can again be seen as a
discretization of the calculations for the relations of a BV
algebra up to homotopy on the chains of the operad $\Arc$ of
\cite{KLP}.
 In \cite{cact} it is proven, that the operad of framed little
discs is  equivalent to the operad of cacti. Moreover, we gave a
description of  cacti in terms of a bi--crossed product of
spineless cacti and an operad built on the monoid $S^1$ which we
showed to be homotopy equivalent to the semi--direct product of
these operads \cite{cact}. Furthermore, we gave a chain model for
spineless cacti in terms of normalized spineless cacti which we
showed to give a natural solution to Deligne's conjecture
\cite{del}. Using the description in terms of the bi--crossed and
semi--direct products, we obtain a chain model for the operad of
framed little discs, by tensoring the chains of normalized
spineless cacti with the chains for the operad built on the monoid
$S^1$. In order to prove the necessary relations on the chain
level one can translate the respective relations from the
relations in the $\Arc$ operad using the method described in
\cite{cact,KLP}. As it turns out, in order to translate the
relations and thus to establish the homotopy BV structure on the
chain level, one needs a refinement of the cell decomposition on
the semi-direct product to be able to accommodate all the
operations which were used in the $\Arc$ operad picture. This
refinement uses cell decompositions on the $S^1$ factors
which are induced by regarding them as the lobe they represent. This leads
to a combinatorial description in terms of planar planted black
and white (b/w) bipartite trees with additional data called
spines. In the language of cacti \cite{cact}, the additional data
keeps track of the position of the local zeros. On these trees,
there are linear orders at each vertex, which may differ from the
induced linear order of the planar planted trees. This forces us
to look at non--rooted trees or equivalently to invert the
orientation of edges. According to the general calculus for
``correlation functions'' defined by trees, to achieve such an
inversion one needs to have a non--degenerate pairing, which is
symmetric and invariant. This is the assumption we have to make on
our algebra. With this assumption, we can rewrite the action of
the cellular chains as ``operadic correlation functions'' for
decorated trees. In this description the operation of the chains
of the framed little discs operad becomes apparent.

The results and techniques we present below can also be employed
in other situations, which we comment on at the end of the paper.
Notably one can use it to obtain an action of cells of a ribbon
graph cell decomposition of moduli space on cyclic complexes. This
should ultimately lead to string topology like operations of the
cells of moduli space of decorated bordered surfaces on the free
loop space of a compact manifold extending the operations of the
string PROP or dioperad. The basic constructions for this are
announced below.

\section*{Acknowledgments}
We would like to thank Alain Connes for an enlightening discussion
and Jim Stasheff for his valuable comments. We also thank the
Max--Planck--Institute for Mathematics in Bonn for providing the
atmosphere and stimulus to conceptualize and complete this paper.

\section{Background}

\subsection{Graphs}
\label{Graphs}

In this section, we formally introduce the graphs and the
operations on graphs which we will use in our analysis of cacti.
This is the approach as given in Appendix B of \cite{cact} in
which cacti are characterized as a certain type of ribbon graph.
Namely, a cactus is a marked treelike ribbon graph with a metric.

\subsubsection{Graphs} A graph $\Gamma$ is a tuple
$(V_{\Gamma},F_{\Gamma}, \imath_{\Gamma}: F_{\Gamma}\rightarrow
F_{\Gamma},\del_{\Gamma}:F_{\Gamma} \rightarrow V_{\Gamma})$ where
$\imath_{\Gamma}$ is an involution $\imath_{\Gamma}^2=id$ without
fixed points. We call $V_{\Gamma}$ the vertices of $\Gamma$ and
$F_{\Gamma}$ the flags of $\Gamma$. The edges $E_{\Gamma}$ of
$\Gamma$ are the orbits of the flags under the involution
$\imath_{\Gamma}$. A directed edge is an edge together with an
order of the two flags which define it.
In case there is no risk of confusion, we will drop the subscripts
$\Gamma$.
Notice that $f\mapsto (f,\imath(f))$ gives a bijection between
flags and directed edges.

We also call $F_v(\Gamma):=\del^{-1}(v)\subset F_{\Gamma}$ the set
of flags of the vertex $v$ and call $|F_v({\Gamma})|$ the valence
of $v$ and denote it by $\val(v)$. We also let
$E(v)=\{\{f,\imath(f)\}|f\in F_{v}\}$ and call these edges the edges
incident to $v$.

The geometric realization of a graph is given by considering  each
flag as a half-edge and gluing the half-edges together using the
involution $\imath$. This yields a one-dimensional CW complex
whose realization we call the realization of the graph.

\subsubsection{Trees}
 A graph is connected if its realization is. A graph is
a tree if it is connected and its realization is contractible.

A rooted tree is a pair $(\t,v_0)$ where $\t$ is a tree and
$v_0\in V_{\t}$ is a distinguished vertex. In a rooted tree there
is a natural orientation for edges, in which the edge points
toward the root. That is we say $(f,\imath (f))$ is naturally
oriented if $\del(\imath(f))$ is on the unique shortest path from
$\del(f)$ to the root. This means that the set $E(v)$ splits up
into incoming and outgoing edges. Given a vertex $v$, we let $|v|$
be the number of incoming edges and call it the arity of $v$. A
vertex $v$ is called a leaf if $|v|=0$. Notice that the root is
the only vertex for which $|v_0|=\val(v_0)$. For all other
vertices $v\neq v_0$ one has $|v|=\val(v)-1$.

A bi-colored or black and white (b/w) tree is a tree $\t$ together with
a map $\color:V\rightarrow \mathbb{Z}/2\mathbb{Z}$. Such a tree is
called bipartite if for all $f\in
F_{\t}:\color(\del(f))+\color(\del(\imath(f)))=1$, that is edges
are only between black and white vertices. We call the set
$V_w:=\color^{-1}(1)$ the white vertices. If $(f,\imath (f))$ is a
naturally oriented edge, we call the edge white if
$\del(\imath(f))\in V_w$ and denote the set of white edges by
$E_w$. Likewise we call $V_b:=\color^{-1}(0)$ the black vertices
and let $E_b$ be the set of black edges, where a naturally
oriented edge $(f,\imath (f))$ is called black if
$\del(\imath(f))\in V_b$.

The black leaves in a rooted black and white tree are called
tails. The edges incident to the tails are called tail edges and
are denoted $E_{tail}$. For tails, we will only consider those
flags of the tail edges which are not incident to the tail
vertices and call them $F_{tail}$.

\subsubsection{Planar trees and Ribbon graphs}

A ribbon graph is a connected graph whose vertices are of valence
at least two together with a cyclic order of the set of flags of
the vertex $v$ for every vertex $v$.

A  graph with a cyclic order of the flags at each vertex gives
rise to bijections $N_v:F_v\rightarrow F_v$ where $N_v(f)$ is the
next flag in the cyclic order. Since $F=\amalg F_v$ one obtains a map
$N:F\rightarrow F$.
The orbits of the map $N \circ \imath$ are called the cycles or
the boundaries of the graph. These sets have the induced cyclic
order.

Notice that each boundary can be seen as a cyclic sequence of
directed edges. The directions are as follows. Start with any flag
$f$ in the orbit. In the geometric realization go along this
half-edge starting from the vertex $\del(f)$, continue along the
second half-edge $\imath(f)$ until you reach the vertex
$\del(\imath(f))$ then continue starting along the flag
$N(\imath(f))$ and repeat.

A tree with a cyclic order of the flags at each vertex is called
planar. A planar tree has only one cycle $c_0$.

\subsection{Planar planted trees}

A planted planar tree is a rooted planar tree $(\t,v_0)$ together
with a linear order of the set of flags at $v_0$. Such a tree has
a linear order of all flags as follows: Let $f$ be the smallest
element of $\del^{-1}(v_0)$, then every flag appears in $c_0$ and
defining the flag $f$ to be the smallest gives a linear order on
the set of all flags. This linear order induces a linear order  on
all oriented edges and on all un-oriented edges, by restricting to
the edges in the orientation opposite the natural orientation
i.e.\ pointing away from the root. We denote the latter by $\prec$
and its restriction to $E(v)$ or $F(v)$ by $\prec_v$.

We will equivalently consider planar planted trees as defined
above or as a rooted planar trees whose root vertex has valence
one. The bijection in one direction is given by adding a new root
vertex and one new edge such that the induced linear structure on
the old root is the given one. This tree is called
the realization of the planar planted tree. In the other direction the
bijection is
simply given by contracting the unique edge incident to the root,
but retaining the linear order.
In the realization of a planar planted tree, we call the unique
edge incident to the (new) root $v_{root}$
the root edge and denote it by $e_{root}$
and set  $f_{root}$ to be the flag of the root edge which is not incident to the root. Also $E_{root}=\{e_{root}\}, F_{root}=\{f_{root}\}$.

An angle at a vertex $v$ in a planar tree is a pair of two flags
incident to $v$ of which one is the immediate successor of the
other in the cyclic order of $F_v$. There is a bijection between
angles, flags and edges by associating to an angle its bigger flag
and to the latter the unique edge defined by it.

\subsection{The genus of a ribbon graph and its surface}
The genus $g(\Gamma)$ of a ribbon graph $\Gamma$ is given by
$2g(\Gamma)+2=|V_\Gamma|-|E_{\Gamma}|+\#cycles$.

The surface $\Sigma(\Gamma)$ of a ribbon graph $\Gamma$ is the
surface obtained from the realization of $\Gamma$ by thickening
the edges to ribbons. I.e.\ replace each 0-simplex $v$ by a closed
oriented disc $D(v)$ and each 1-simplex $e$ by $e\times I$
oriented in the standard fashion. Now glue the boundaries of
$e\times I$ to the appropriate discs in their cyclic order
according to the orientations. Notice that the genus of
$\Sigma(\Gamma)$ is $g(\Gamma)$ and that $\Gamma$ is naturally
embedded as the spine of this surface.

\subsubsection{Treelike and marked ribbon graphs}

A ribbon graph together with a distinguished cycle $c_0$ is called
{\em treelike} if
\begin{itemize}
\item[i)]
the graph is of genus $0$  and
\item[ii)] for all
cycles $c_i\neq c_0$: if $f\in c_i$ then $\imath(f)\in c_0$.
\end{itemize}
In other words each edge is traversed by the cycle $c_0$.
Therefore there is a cyclic order on all (non-directed) edges,
namely the cyclic order of $c_0$.

 A {\em marked ribbon graph} is a
ribbon graph together with a map $\mk:\{cycles\} \rightarrow
F_{\Gamma}$ satisfying the conditions
\begin{itemize}
\item[i)] For every cycle $c$ the directed edge $\mk(c)$ belongs
to the cycle.

\item[ii)] All vertices of valence two are in the image of $\mk$,
that is $\forall v,\val(v)=2$ implies  $v\in Im(\del\circ\mk)$.
\end{itemize}

Notice that on a marked treelike ribbon graph there is a linear
order on each of the cycles $c_i$. This order is defined by
upgrading the cyclic order to the linear order $\prec_i$ in which
$\mk(c_i)$ is the smallest element.

\subsubsection{Dual b/w tree of a marked ribbon graph}
Given a marked treelike ribbon graph $\G$, we define its dual tree
to be the colored graph whose black vertices are given by $V_{\G}$
and whose set of white vertices is the set of cycles $c_i$ of
$\G$. The set of flags at $c_i$ are the flags $f$ with $f\in c_i$
and the set of flags at $v$ are the flags $\{f:f \in c_0,
\del(f)=v\}$. The involution is given by $\imath_{\t}(f)=N(f)$ if
$f\in c_0$ and $\imath_{\t}(f)=N^{-1}(f)$ else.

This graph is a tree and is b/w and bipartite by construction. It
is also planar, since the $c_i$ and the sets $F(v)$ have a cyclic
order and therefore also $F_v\cap c_0$. It is furthermore rooted
by declaring $\del(\mk(c_0))$ to be the root vertex and declaring
$\mk(c_0)$ to be the smallest element makes it into a planted
tree.

An equivalent definition is given by defining that there is an
edge between a pair of a black and a white vertex if and only if
the vertex corresponding to $b$ is on the boundary of the cycle
$c_i$, i.e.\ $v\in \del(c_i):= \{\del(f):f\in c_i\}$.

\subsubsection{Spineless marked ribbon graphs}
\label{spinlessgraph} A marked treelike ribbon graph is called
{\em spineless}, if

\begin{itemize}
\item[i)] There is at most one vertex of valence $2$. If there is
such a vertex $v_0$ then $\del(mk(c_0))=v_{0}$.

\item[ii)] The induced linear orders on the $c_i$ are compatible
with that of $c_0$, i.e.\ $f\prec_i f'$ if and only if
$\imath(f')\prec_0 \imath(f)$.
\end{itemize}

\subsubsection{Graphs with a metric}

A metric $w_{\Gamma}$ for a graph is a map $E_{\Gamma}\rightarrow
\mathbb{R}_{>0}$. The (global) re-scaling of a metric $w$ by
$\lambda$ is the metric $ \lambda w: \lambda w(e)=\lambda w(e)$.
The length of a cycle $c$ is the sum of the lengths of its edges
$length(c)=\sum_{f\in c} w(\{f,\imath(f)\})$. A metric for a
treelike ribbon graph is called normalized if the length of each
non-distinguished cycle is $1$.

\subsubsection{Marked ribbon graphs  with metric and maps of circles.}
For a marked ribbon graph with a metric, let $c_i$ be its cycles,
let $|c_i|$ be their image in the realization and let $r_i$ be the
length of $c_i$. Then there are natural maps
$\phi_i:S^1\rightarrow |c_i|$ which map $S^1$ onto the cycle by
starting at the vertex $v_i:=\del(\mk(c_i))$ and going around the
cycle mapping each point $\theta\in S^1$ to the point at distance
$\frac{\theta}{2\pi}r_i$ from $v_i$ along the cycle $c_i$.

\subsubsection{Contracting edges}
The contraction $(\bar V_{\Gamma}, \bar F_{\Gamma},\bar
\imath,\bar \del)$ of a graph
$(V_{\Gamma},F_{\Gamma},\imath,\del)$ with respect to an edge
$e=\{f,\imath(f)\}$ is defined as follows. Let $\sim$ be the
equivalence relation induced by $\del(f)\sim\del(\imath(f))$. Then
let $\bar V_{\Gamma}:=V_{\Gamma}/\sim$, $\bar
F_{\Gamma}=F_{\Gamma}\setminus\{f,\imath(f)\}$ and $\bar \imath:
\bar F_{\Gamma}\rightarrow \bar F_{\Gamma}, \bar\del: \bar
F_{\Gamma}\rightarrow \bar V_{\Gamma}$  be the induced maps.

For a marked ribbon graph, we define the marking of $(\bar
V_{\Gamma}, \bar F_{\Gamma},\bar \imath,\bar \del)$ to be
$\overline{\mk}(\bar c)=\overline{\mk(c)}$ if
$\mk(c)\notin\{f,\imath(f)\}$ and $\overline{\mk}(\bar
c)=\overline{N\circ \imath(\mk (c))}$ if
$\mk(c)\in\{f,\imath(f)\}$, viz.\ the image of the next flag in
the cycle.

\subsubsection{Labelling graphs}
By a labelling of the edges of a graph $\Gamma$ by a set $S$, we
simply mean a map $E_{\Gamma}\rightarrow S$. A labelling of a
ribbon graph $\Gamma$ by a set $S$ is a map $\lab\{$cycles of
$\Gamma\}\rightarrow S$, we will write $c_i:=\lab^{-1}(i)$. By a
labelling of a black and white tree  by a set $S$ we mean  a map
$\lab:E_w\rightarrow S$. Again we will write $v_i:=\lab^{-1}(i)$.

\subsubsection{Planar planted bipartite labelled trees with white leaves}
We set $\wlbptree(n)$ to be the set of planar planted bipartite
trees which are labelled from $\{1,\dots,n\}$ with  white leaves
only. To avoid cluttered notation, we
 also denote the respective free Abelian group and the $k$-vector
space with basis $\wlbptree(n)$ by the same name and let
 $\wlbptree$ be their union respectively direct sum.

 \subsection{Cacti}

\begin{df}
A cactus with $n$ lobes is  a $\{0,1, \dots ,n\}$ labelled marked
treelike ribbon graph with a metric. I.e.\ The set  $\Cacti(n)$ is
the set of these graphs. $\Cact(n)\subset \Cacti(n)$ is the subset
of spineless graphs and called the spineless cacti or
alternatively cacti without spines. $\Cacti^1(n)\subset \Cacti(n)$
is the subset of normalized graphs, called normalized cacti, and
finally $\Cact^1(n)=\Cact(n)\cap\Cacti^1(n)$ is the set of
normalized spineless cacti.
\end{df}

\subsubsection{Cactus terminology}
The edges of a cactus are traditionally called arcs or segments
and the cycles of a cactus are traditionally called lobes. The
vertices are sometimes called the marked or special points.
Furthermore the distinguished cycle $c_0$ is called the outside
circle or the perimeter and the vertex $\del(\mk(c_0))$ is called
the global zero. And the vertices $\del(\mk(c_i)),i\neq 0$ are
called the local zeros. In pictures these are represented by lines
rather than fat dots.

\begin{rmk}
\label{setrem} It is clear that as sets $\Cacti(n)=\Cact(n)\times
(S^1)^{\times n}$ and $\cact(n)= \cact^1(n)\times
\mathbb{R}_{>0}^{\times n}$.

For the first statement one notices for each lobe $v_i$ there is a
unique lowest intersection point $b$ which is the vertex of the
outgoing edge of $v$. Thus there is a canonical map
$\phi'_i:S^1\rightarrow |c_i|$ which starts at $b$ and goes around
the cycle opposite its natural orientation. So to each cycle we associate
$(\phi'_i)^{-1}(\del(\mk(c_i)))$ that is the co-ordinate of the
spine as measured by $\phi'_i$. This gives the projection onto the
factors $(S^1)^{\times n}$. The projection onto the first factor
is given by forgetting the spines, i.e.\ contracting the edges
$\mk(c_i)$ if $\val(\del(\mk(c_i)))=2$ and changing the marking
 to the unique marking which makes the graph spineless.

 For the second statement the first projection is given by
 homogeneously scaling the weights of the edges of each non-marked
 cycle so that their lengths are one. The projection to the
 factors of $\mathbb{R}_{>0}$ are given by associating to each
 lobe its length.
 In both cases the inverse map is clear.
\end{rmk}

\begin{df}
The topological type of a spineless cactus in $\cact^1(n)$ is
defined to be its dual b/w tree $\t \in \wlbptree(n)$.
\end{df}

\begin{rmk}
\label{arctoedge}
 Notice that the arcs of a  cactus correspond to
the set  $E_{arcs}=E(\t)\setminus (\{e_{root}\})$. This bijection
can be defined as follows. To a given $e\in E_{arcs}, e=\{w,b\}$
with $b$ black and $w$ white, we associate the unique arc
 between the points corresponding to the
black vertices $b$ and $b-$ where $b-$ is the black vertex
immediately preceding $b$ in the cyclic order of $v$. In other
words if $e=\{f,\imath(f)\}$ with $f\in F_v$. Let $f-$ be the flag
immediately preceding $f$ in the cyclic order at $v$, then
$b-=\del(\imath(f-))$. Notice that if $|v|=0$ then and only then
$f-=f$.
\end{rmk}

\begin{rmk}
\label{typelemma} A spineless cactus is uniquely determined by its
topological type and the lengths of the segments.
\end{rmk}

\subsection{The CW complex of normalized spineless cacti}
We recall from \cite{del} the CW complexes $\CWcact(n)$. For more details
and pictures the reader is referred to \cite{del,cact}.

\begin{rmk}

\label{lengthrem} For a normalized spineless cactus the lengths of
the arcs have to sum up to the radius of the lobe and the number
of arcs on a given lobe  represented by a white vertex $v$ is
$\val(v)=|v|+1$. Hence the lengths of the arcs lying on the lobe
represented by a vertex $v$ are in 1-1 correspondence with points
of the simplex $|\Delta^{|v|}|$. The coordinates of
$|\Delta^{|v|}|$ naturally correspond to the arcs of the lobe
represented by $v$ on one hand and on the other hand in the dual b/w graph
to the edges incident to $v$.
\end{rmk}

\subsubsection{The tree differential in the spineless case}
\label{diffdef}
 Let $\t\in \wlbptree$. We set $E_{angle}=E(\t)\setminus
(E_{leaf}(\t)\cup \{e_{root}\})$ and we denote by
$\num_E:E_{angle} \rightarrow \{1,\dots,N\}$  the bijection which
is induced by the linear order $\prec^{(\t,p)}$.

 Let $\t\in \wlbptree$, $e\in E_{angle}$,
$e=\{w,b\}$, with $w\in V_w$ and $b\in V_b$. Let $e-=\{w,b-\}$ be
the edge preceding $e$ in the cyclic order $\prec^{\t}_w$ at $w$.
Then $\del_e(\tau)$ is defined to be the planar tree obtained by
collapsing the angle between the edge $e$ and its predecessor in
the cyclic order of $w$ by identifying $b$ with $b-$ and $e$ with
$e-$.
 Formally
$w=\whitevert(e), e-=\prec^{\t}_w(e),\{b-\}= \del(e-)\cap
V_b(\t)$, $V_{\del_e(\tau)}=V(\t)/(b\sim b-)$,
$E_{\del_e(\tau)}=E_{\tau}/(e\sim e-)$.
The linear order of $\del_e(\t)$ is given by keeping the linear
order at all vertices which are not equal to $\bar b$ where $\bar
b$ is the image of $b$ and $b-$. For $\bar b$ the order is given
by extending the linear order $(\In(\bar b), \prec_{\bar
b}^{\del_e(\t)}) =(\In(b-)\amalg\In(b), \prec^{\t}_{b-}\amalg
\prec^{\t}_{b}) $ ---the usual order on the union of totally
ordered sets-- to $E(\bar b)$ by declaring the image of $e$ and
$e-$ to be the minimal element.

\begin{df}
We define the operator $\del$ on the space $\wlbptree$ to be given
by the following formula: $\del(\t) := \sum_{e\in E_{angle}}
(-1)^{\num_E(e)-1} \del_e (\tau) $.
\end{df}

\subsubsection{The Cell Complex}

\begin{df}
We define $\wlbptree(n)^k$ to be the elements of $\wlbptree(n)$
with $|E_w|=k$.
\end{df}

\begin{df} For $\t \in \wlbptree$ we define
$\D(\t):=\times_{v \in V_w(\tau)}\D^{|v|}$. We define
$C(\t)=|\D(\t)|$.  Notice that $\dim(C(\t))=|E_w(\t)|$.

Given $\D(\t)$ and a vertex $x$ of any of the constituting
simplices of $\D(\t)$ we define the $x$-th face of $C(\t)$ to be
the subset of $|\D(\t)|$ whose points have the $x$-th coordinate
equal to zero.
\end{df}

\begin{df}
We let $\CWcact(n)$ be the CW complex whose k-cells are indexed by
$\t \in \wlbptree(n)^k$ with the cell $C(\t)=|\D(\t)|$ and the
attaching maps $e_{\t}$ defined as follows. We identify the $x$-th
face of $C(\t)$ with $C(\t')$ where $\t'=\del_x(\t)$. This
corresponds to contracting an edge of the cactus if its weight
goes to zero (see Remark \ref{arctoedge}) so that $\Delta(\del \t)$
is identified with $\del (\Delta(\t))$.
\end{df}

\begin{df}
We define the topology of $\cact^1(n)$ to be that induced by the
bijection with $\CWcact(n)$. Via Remark \ref{setrem} this gives a
topology to the spaces $\Cact(n),\cacti(n)$ and $\cacti^1(n)$.
\end{df}

\subsection{The (quasi)-operad structure}
\subsubsection{The operad of cacti}
The gluing maps for cacti
\begin{equation}
\circ_i:\cacti(n)\otimes \cacti(m)\rightarrow \cacti(n+m-1)
\end{equation}
are defined on elements $(c,c')\mapsto c\circ_i c'$ as follows
\begin{itemize}
\item[1)] Scaling the weight function $w'$ of $c'$ by the length
$\frac{r_i}{R}$  where $r_i$ is the length of the cycle $c_i$ of
the cactus $c$ and $R$ is the length of the cycle $c_0$ of $c'$.
\item[2)] Identifying the realization of the cycle $c_0$ of $c'$
with the cycle $c_i$ of $c$ via the maps $\phi_0(c')$ and
$\phi_i(c)$, with the orientation on the second $S^1$ reversed, as
usual.
\end{itemize}

These maps together with the $\Sn$ action permuting the labels
turn the collection $\{\cacti(n)\}$ into an operad  $\cacti$. The
collection $\{\cact(n)\}$ forms the suboperad $\cact$.

\subsubsection{The quasi-operad of normalized cacti}
We recall from \cite{cact} that a quasi-operad is the
generalization of a (pseudo)-operad in which the axiom of
associativity is omitted and the others are kept.

The gluing maps for normalized cacti
\begin{equation}
\circ_i:\cacti^1(n)\otimes \cacti^1(m)\rightarrow \cacti^1(n+m-1)
\end{equation}
are defined on elements $(c,c') \mapsto c\circ_i c'$ simply by
identifying the realization of the cycle $c_0$ of $c'$ with the
cycle $c_i$ of $c$ via the maps $\phi_0(c')$ and $\phi_i(c)$ again
with the orientation on the second $S^1$ reversed.

These maps together with the $\Sn$ action permuting the labels
turn the collection $\{\cacti^1(n)\}$ into a homotopy associative
quasi-operad  $\cacti^1$. The collection $\{\cact^1(n)\}$ forms a
homotopy associative quasi-suboperad $\cact^1$ of $\cacti^1$ \cite{cact}.

\subsection{Relations among cacti}
\begin{thm}
\label{cactthm}
\cite{cact}
Normalized cacti are homotopy equivalent through quasi-operads to
the cacti. The same holds for the (quasi)-suboperads of normalized spineless cacti and
spineless cacti.
\end{thm}

\begin{cor}\cite{cact}
Normalized cacti are quasi-isomorphic as quasi-operads to cacti and
normalized spineless cacti are quasi-isomorphic as quasi-operads
to spineless cacti. In particular in both cases the homology quasi-operads
are operads and are isomorphic as operads.
\end{cor}

\subsubsection{Remarks on the bi-crossed product}
In this section we recall the construction of the bi-crossed
product as it was given in \cite{cact} to which we refer the
reader for more details.

First notice that there is an action of $S^1$ on $\Cact(n)$ given
by rotating the base point {\em clockwise} (i.e.\ in the
orientation opposite the usual one of $c_0$)  around the
perimeter. We denote this action by
$$\rho^{S^1}: S^1 \times \Cact(n) \rightarrow \Cact(n)$$
With this action we can define the twisted gluing
\begin{eqnarray}
\label{circtheta}
\circ_i^{S^1}:\Cact(n) \times S^1(n) \times \Cact(m) &\rightarrow& \Cact(n+m-1)\nn\\
(C,\theta,C')&\mapsto& C \circ \rho^{S^1}(\theta_i,C') =: C
\circ_i^{\theta_i}C'
\end{eqnarray}

Given a cactus  without spines $C\in \Cact(n)$ the orientation
reversed perimeter (i.e.\ going around the outer circle {\em
clockwise} i.e.\ reversing the orientation of the source of
$\phi_0$) gives a map $\Delta_C: S^1 \rightarrow (S^1)^n$.

As one goes around the perimeter the map goes around each circle
once and thus the map $\Delta_C$ is homotopic to the diagonal
$ \Delta_C (S^1) \sim \Delta(S^1)$.

We can use the map $\Delta_C$ to give an action of $S^1$ and
$(S^1)^{\times n}$.
\begin{equation}
\rho^C: S^1 \times(S^1)^{\times n}\stackrel{\Delta_C}
{\rightarrow} (S^1)^{\times n} \times (S^1)^{\times n}
\stackrel{\mu^n}{\rightarrow}(S^1)^{\times n}
\end{equation}
here $\mu_n$ is the diagonal multiplication in $(S^1)^{\times n}$
and $\bar \circ_i$ is the operation which forgets the $i$-th
factor and shuffles the last $m$ factors to the $i$-th, \dots ,
$i+m-1$st places. Set
\begin{multline}
\label{perturbdef} \circ_i^C:(S^1)^{\times n} \times (S^1)^{\times
m} \stackrel{(id \times \pi_i)(\Delta) \times id}
{\longrightarrow} (S^1)^{\times n} \times
S^1\times (S^1)^{\times m}\\
\stackrel{id \times \rho^C}{\longrightarrow} (S^1)^{\times n}
\times (S^1)^{\times m}
\stackrel{\bar\circ_i}{\longrightarrow}(S^1)^{\times n+m-1}
\end{multline}
These maps are to be understood as perturbations of the usual maps
\begin{multline}
 \circ_i:(S^1)^{\times n} \times (S^1)^{\times
m} \stackrel{(id \times \pi_i)(\Delta) \times id}
{\longrightarrow} (S^1)^{\times n} \times
S^1\times (S^1)^{\times m}\\
\stackrel{id \times \rho}{\longrightarrow} (S^1)^{\times n} \times
(S^1)^{\times m}
\stackrel{\bar\circ_i}{\longrightarrow}(S^1)^{\times n+m-1}
\end{multline}
where now $\rho$ is the diagonal action of $S^1$ on $(S^1)^{\times
n}$. The maps $\circ_i$ and the permutation action on the factors
give the collection $\{\mathcal{S}^1(n)\}=(S^1)^{\times n}$ the structure
of an operad. In fact this is exactly the usual construction of an
operad built on a monoid.

\begin{thm}
\label{cactbicross} \cite{cact}  The operad of cacti is the
bi--crossed product of the operad $\cact$ of  spineless cacti with
the operad $\mathcal {S}^1$ based on $S^1$. Furthermore this
bi--crossed product is homotopic to the semi--direct product of
the operad of cacti without spines with the circle group $S^1$.
\begin{equation}
\cacti \cong \cact \bowtie {\mathcal S}^1 \simeq \cact \rtimes
{\mathcal S}^1
\end{equation}
The multiplication in the bi-crossed product is given by
\begin{equation}
(C,\theta) \circ_i (C',\theta') = (C\circ_i^{\theta_i} C',
\theta\circ_{i}^{C'}\theta')
\end{equation}
The multiplication in the semi-direct product is given by
\begin{equation}
(C,\theta) \circ_i (C',\theta') = (C\circ_i^{\theta_i} C',
\theta\circ_{i}\theta')
\end{equation}
Also, normalized cacti are homotopy equivalent to cacti which are
homotopy equivalent to the bi-crossed product of normalized cacti
with $\mathcal{S}^1$ and the semi-direct product with
$\mathcal{S}^1$, where all equivalences are as quasi-operads
\begin{equation}
\cacti^1 \sim \cacti \cong \cact \bowtie {\mathcal S}^1
\sim\cact^1 \bowtie {\mathcal S}^1\sim  \cact^1 \rtimes {\mathcal
S}^1
\end{equation}
\end{thm}

\begin{rmk}
The proof of the first statement is given by verifying that the
two operad structures coincide. For the second statement one
notices that the homotopy diagonal is homotopy equivalent to the
usual one and that one can find homotopies to the diagonal which
continuously depend on the cactus. The third statement follows
from contracting the factors $\mathbb{R}^n_{>0}$ and using
Theorem \ref{cactthm}.
\end{rmk}

\begin{cor}
The homology operad of $\cacti$ is the semi-direct product of
$\cacti$ and the homology of the operad $\mathcal{S}^1$ built on
the monoid $S^1$.
\end{cor}

\subsection{Relation to (framed) little discs}

\begin{thm}\cite{cact}
The operad $\cact$ is equivalent to the little discs operad and
the operad  $\cacti$ is equivalent to the framed little discs
operad.
\end{thm}

The latter result has been first stated by Voronov in \cite{Vor}.

\section{A CW decomposition for $\cacti^1$ and a chain model for the framed little discs}

\begin{df}
A $\Zz$ decoration for a black and white bipartite tree is a map
$\Zdec: V_w \rightarrow \Zz$.
\end{df}

\begin{prop}
\label{firstcells} The quasi--operad of normalized cacti
$\cacti^1$ has a CW--decomposition which is given by cells indexed
by planar planted bi--partite trees with a  $\Zz$ decoration. The
$k$ cells are indexed by trees with $k-i$ white edges and $i$
vertices marked by $1$.

Moreover cellular chains are a chain model for the framed little
discs operad and form an operad. This operad is isomorphic to the
semi--direct product of the chain model of the little discs operad
given by $CC_*(\cact)$ of \cite{del} and the cellular chains of
the operad built on the monoid $S^1$.
\end{prop}

\begin{proof}
For the CW decomposition we note that  as spaces $\cacti^1(n)=
\cact^1(n) \times (S^1)^{\times n}$ see Remark \ref{setrem}. Now
viewing $S^1=[0,1]/0\sim1$ as a 1-cell together with the 0-cell
given by $0\in S^1$ the first part of the proposition follows
immediately, by viewing the decoration by 1 as indicating the
presence of the 1-cell of $S^1$ for that labelled component in the
product of cells.

To show that the cellular chains indeed form an operad, we use the
fact that the bi--crossed product is homotopy equivalent to the
semi--direct product in such a way, that the action of a cell
$S^1$ in the bi--crossed product is homotopic to the diagonal
action. This is just the observation that the diagonal and the
diagonal defined by a cactus are homotopic. Since a semi-direct
product of a monoid with an operad is an operad the statement
follows. Alternatively one could just remark, that there is also
an obvious functorial map induced by the diagonal for these cells.

The chains are a chain model for the framed little discs operad
since  $\cacti^1(n)$ and $\cacti(n)$ are homotopy equivalent and
the latter is equivalent to the framed little discs operad.
\end{proof}

Although the above chain model is the one one would expect to use
for framed little discs, it does not have enough cells for our
purposes. In order to translate the proofs in the arc complex
given in \cite{KLP} into statements about the Hochschild complex,
we will need a slightly finer cell structure then the one above.
After having used the larger structure one can reduce to the cell
model with less cells as they are obviously equivalent.

\begin{figure}
\epsfxsize = \textwidth
\epsfbox{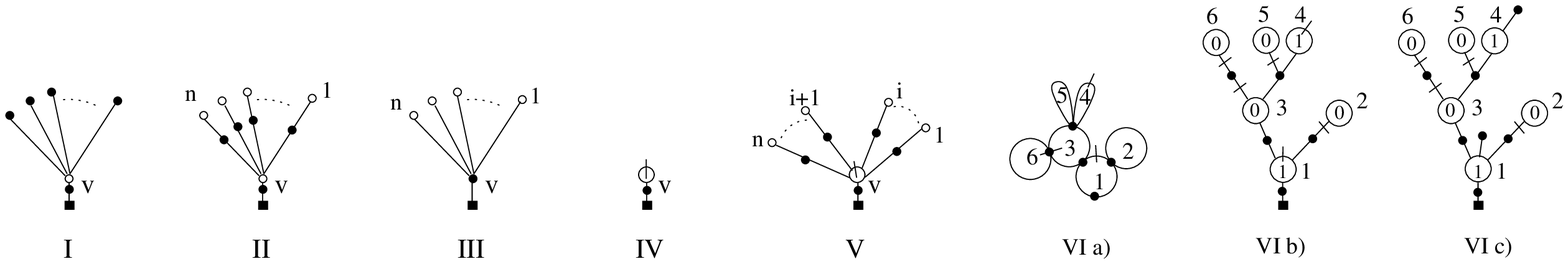}
\caption{\label{cactexamples}I. the tree $l_n$ II. the tree $\t_n$
III. the tree $\t_n^b$ IV. the tree O'.V. the tree $\t'_{n,i}$ VI.a) a marked
tree-like ribbon graph b) the corresponding decorated tree
c) its realization}
\end{figure}

\begin{df}
A spine decoration $\sdec$ for a planted planar bi--partite tree
is a $\Zz$ decoration together with the marking of one angle at
each vertex labelled by one and a flag at each vertex labelled by
zero.  We call the set of such trees which are $n$-labelled by
$\swlbptree(n)$ and again use this notation as well for the free
Abelian group and the $k$ vector space generated by these sets. We
let
 $\swlbptree$ be their union respectively direct sum.
In pictures we show the angle marking as a line emanating from the
vertex which lies between the marked edges and an edge marking by
a line through the respective edge.  For an example see Figure
\ref{cactexamples} VI. We sometimes omit the edge marking if the
marked edge is the outgoing edge, e.g.\ in Figure
\ref{bvtopartcact}.

A realization $\hat \t$ of a planar planted bi--partite tree $\t$
with a spine decoration is a realization of $\t$ as a planar
planted tree (the root is fixed to be black) together with one
additional edge inserted into each marked angle connecting to a
new vertex. We call the set of these edges spine edges and denote
them by $E_{spine}$. Likewise set $V_{spine}$ to be the set of new
vertices called the spine vertices which are defined to be black.
The spine edges are then white edges. Like for
tails, we will only consider the flags of $E_{spine}$, which are
not incident to the spine vertices. We call the set of these flags
$F_{spine}$. Notice that this tree is the dual tree of a cactus
with an explicit marking of the flags $\mk(c_i)$. Given a cactus,
we call its dual tree with explicit markings its topological type.
If $\t$ had tails, we will split the set of tails of the
realization into spines and free tails which are the images of the
original tails. $E_{tails}(\hat\t)=E_{ftails}(\hat \t)\amalg
E_{spine}(\hat \t)$ and likewise for the respective flags.

A spine decoration induces a new linear order on the flags
incident to the white vertices of its realization. This order
$\prec'_v$ is given by the cyclic order at $v$ and declaring the
smallest element to be the spine flag in case $\Zdec(v)=1$ and the
marked flag in case $\Zdec(v)=0$. This gives a canonical
identification of $\Fnum:F_v \rightarrow \{0,\dots, |v|\}$.
\end{df}

\begin{prop}
\label{secondcells} The spaces $\cacti^1(n)$ of the quasi--operad
of normalized cacti $\cacti^1$ have  CW--decompositions
$\CWtwo(n)$ whose  cells are indexed by spine decorated planar
planted bi--partite trees $(\t,\sdec)\in \swlbptree$ corresponding
to the topological type of the cacti. The $k$ cells are indexed by
$n$-labelled trees with $k-i$ white edges and $i$ markings by $1$.

Moreover cellular chains of the complex above are a chain model
for the framed little discs operad and form an operad.
\end{prop}

\begin{proof}
The decomposition is almost as in the preceding proposition except
that in the product $\cact^1(n)\times (S^1)^{\times n}$ we
decompose each factor $S^1$ as indicated by the lobe it presents.
I.e.\ for the $S^1$ associated to the $n$--th lobe we chose the
0--cells to be corresponding to the marked points and 1--cells
corresponding to the arcs with gluing given by attaching the
1--cells to the 0--cells representing the endpoints of the arcs.
(E.g. 4 0-cells and 4 1-cells for the lobe 1 in Figure \ref{cactexamples} VIa).
In terms of trees, the arcs correspond to the angles and thus we
take a marking of an arc to be the inclusion of the corresponding
1-cell in the tensor product of the cell complexes. Likewise the
edges correspond to the marked points and we take a marking of an
edge to be the inclusion of the corresponding 0-cell in the tensor
product of the cell complexes.

For the operadic properties, we remark that moving the spine along
an arc and then gluing, which is what is parameterized by marking
an angle on the lobe $i$ of $c$ when calculating $c\circ_ic'$, has
the effect of moving the base point of $c'$ along a complete
sequence of arcs until it coincides with a marked point in the
composition of the two cacti. This is one side of the bi-crossed
product.  The effect on the local zeros of $c'$ of the movement of
the base point is to move them corresponding to structure maps of
the bi-crossed product above. The local zeros thus move through a
full arc if the global zero passes through the arc on which they
lie.
Therefore the $\circ_i$ product of two cells results in sums of
cells. Marking an arc of $c'$ obviously gives rise
to a sum of cells. Alternatively, one can again just remark that there is a
functorial map for the diagonal for this cell model, since there
is such a map on the first factor by \cite{del} and its existence
is obvious on the second factor.

The associativity follows from the associativity of cacti. Let
$C(\t)$, $\t\in \swlbptree(n)$ be the cells in the CW-complex and
$\dot C(\t)$ their interior. Then $P(\t)=\dot C(\t)\times
\mathbb{R}_{>0}^n, \t\in\swlbptree$ give a pseudo-cell
decomposition $Cacti(n)=\amalg_{\t} P(\t)$.
It is easy to see that $Im(P(\t)\circ_i
P(\t'))=\amalg_k P(\t_k)$ for some $\t_k$
 and $\circ_i$ is a bijection onto its image.
Let $\circ_i^{comb}$ be the quasi-operad structure pulled back
from $\CWtwo$ to $\swlbptree$ and $\circ_i^{+}$ be the operad
structure pulled back from the pseudo-cell decomposition of
$\Cacti$ to $\swlbptree$. Then these two operad structures
coincide over $\mathbb{Z}/2\mathbb{Z}$ thus yielding associativity
up to signs. The signs are just given by shuffles, c.f.\
\S\ref{signsection}, and are associative as well.
\end{proof}

\begin{rmk} Pulling back the operadic compositions, the differential
and the grading yields a dg-operad structure on $\swlbptree$ which
is isomorphic to that of $\CCCacti:=\bigoplus_n CC_*(\CWtwo(n))$.

The operation is  briefly as follows: given two trees
$\t,\t'\in\swlbptree$ the product is $\t\circ^{comb}_i\t'=\sum \pm
\t_k$ where the $ \t_k$ are the trees obtained by the following
procedure. Delete $v_i$ to obtain an ordered collection of trees
$(\t^c_l,\prec'_v)$ then graft these trees to $\t'$ keeping their
order by first identifying the spine edge or marked edge of $v_i$
with the root edge of $\t'$ and then grafting the rest of the
branches to $\t'$ so that their original order is compatible with
that of $\t'$. Lastly contract the image of the root edge of $\t'$
and declare the image of the root of $\t$ to be the new root. The
sign is as explained in \ref{signsection}. Due to the isomorphism
between $\CCCacti$ and $\swlbptree$ we will drop the superscript
$comb$.
\end{rmk}

\subsection{The GBV structure}
The picture for the GBV structure is essentially that of \cite{KLP} and goes back to
\cite{CS1}.
It appears here is another guise, however, since we are now dealing with cells
in $\CCCacti$.

\begin{figure}
\epsfxsize =4in \epsfbox{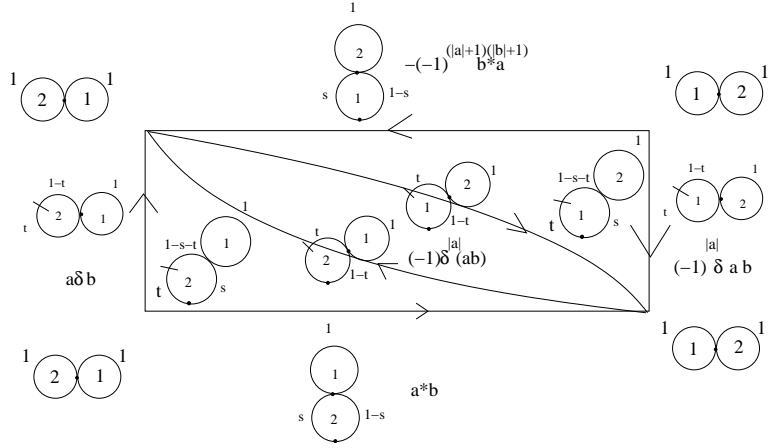} \caption{\label{bvcactfig}
The compatibility of the BV operator and the bracket}
\end{figure}

First notice that there is a product on the chain level induced by
the spineless cactus given by the rooted tree $\t_n$ depicted in
Figure \ref{cactexamples}. Explicitly: $a\cdot b\mapsto
\gamma(\t^b_2; a,b)$ where $\gamma$ is the usual operadic
composition. This product gives $\CCCacti$ the structure of an
associative algebra with unit. Moreover the product is commutative
up to homotopy. The homotopy is given by the usual operation which
is induced by $\gamma(\t_1;a,b)$. This also induces a bracket
which is Gerstenhaber up to homotopy. This can be seen by
translating the statements from \cite{KLP,del}, but it  also
follows from the BV description of the bracket below (Figure
\ref{bvcactfig}).

\begin{figure}
\epsfxsize = 4in \epsfbox{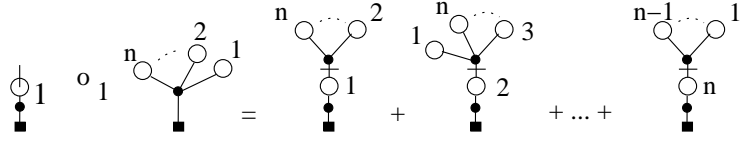}
\caption{\label{bvtopartcact}The decomposition of the BV operator}
\end{figure}

To give the BV structure, let $O'$ be the tree  with one white
vertex, no additional black edges, no free tails and a spine.
Notice that the operation $\delta$ induced by $a\mapsto \gamma(
O',a)$ on $\CCCacti$  breaks up on products of chains as follows,
see Figure \ref{bvtopartcact}
\begin{eqnarray}
\label{threedel}
\delta(ab) &\sim& \vardel(a,b) + (-1)^{|a||b|}\vardel(b,a)\nn\\
\delta(abc) &\sim& \vardel(a,b,c)+(-1)^{|a|(|b|+|c|)}\vardel(b,c,a)\nn\\
&&+(-1)^{|c|(|a|+|b|)}\vardel(c,a,b)
\end{eqnarray}
\begin{eqnarray}
\delta(a_1 a_2\cdots a_n)&\sim& \sum_{i=0}^{n-1}
(-1)^{\s(c^i,a)} \delta(a_{c^i(1)}, \dots, a_{c^i(n)})
\end{eqnarray}
where $c$ is the cyclic permutation and $\s(c^i,a)$ is the sign
of the cyclic permutations of the graded elements $a_i$.

\begin{lem}
\label{partbv}
\begin{equation}
\vardel(a,b,c) \sim (-1)^{(|a|+1)|b|} b \vardel(a,c)
+\vardel(a,b)c -\vardel(a)bc
\end{equation}
\end{lem}

\begin{figure}
\epsfxsize = 4in \epsfbox{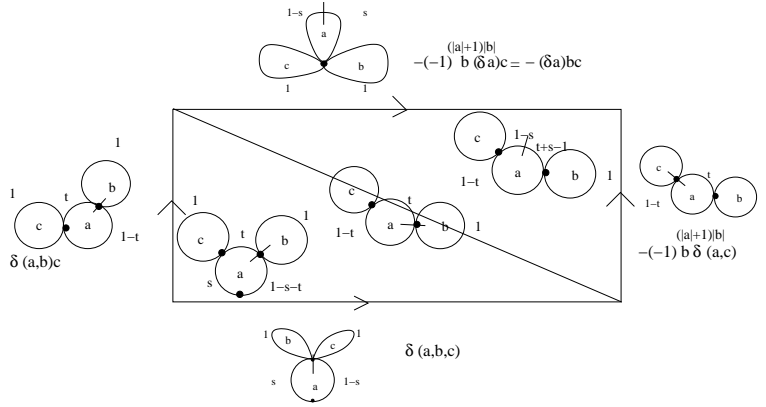}
\caption{\label{cactpartbvfig}The basic chain homotopy responsible
for BV}
\end{figure}

{\bf Proof.} The proof is contained in Figure \ref{cactpartbvfig}.

\begin{prop}
\label{GBVprop}
The chains $\CCCacti$ are a GBV algebra up to homotopy.
\end{prop}

\begin{proof}
The BV structure follows from the Lemma \ref{partbv} via the
calculation:
\begin{eqnarray}
\delta(abc) &\sim& \vardel(a,b,c)+(-1)^{|a|(|b|+|c|)}\vardel(b,c,a)
+(-1)^{|c|(|a|+|b|)}\vardel(c,b,a)\nn\\
&\sim&(-1)^{(|a|+1)|b|} b \vardel(a,c) +\vardel(a,b)c
-\vardel(a)bc + (-1)^{|a|} a \vardel(b,c)\nn
\end{eqnarray}
\begin{eqnarray}
 &&+(-1)^{|a||b|}\vardel(b,a)c -(-1)^{|a|}a\vardel(b)c
+(-1)^{(|a|+|b|)|c|} a \vardel(b,c)\nn\\
&&+(-1)^{|b|(|a|+1|)+|a||c|}b\vardel(c,a)c - (-1)^{|a|+|b|} ab\vardel(c)\nn\\
&\sim&\delta(ab)c+(-1)^{|a|}  a\delta(bc) + (-1)^{|a+1||b|} b\delta (ac)
 -\delta(a)bc\nn\\
&&-(-1)^{|a|} a\delta(b)c-(-1)^{|a|+|b|}ab\vardel(c)
\end{eqnarray}
Figure \ref{bvcactfig} contains the homotopy relating the BV
operator to the bracket.
\end{proof}

\section{The action}

\subsection{Assumption}
\label{assumption}
Now we fix $A$ to be a finite dimensional
 associative algebra with unit $1$ together
with an  inner product $\eta: A\otimes A\rightarrow k$ which is
non-degenerate and both i) invariant: $\eta(ab,c)=\eta(a,bc)$ and
ii) symmetric: $\eta(a,b)=\eta(b,a)$. Such an algebra is called a
Frobenius algebra.

We will use $CH$ to stand for Hochschild cochains $CH^n(A,A):=
Hom(A^{\otimes n},A)$.

Actually, it would be enough to have a non-degenerate
inner-product $\eta$ on $A\simeq CH^0(A,A)$ for which  i) holds on
$HH^0(A,A)$, that is up to homotopy for $A$. The condition ii)
will then hold automatically up to homotopy since $CH^0(A,A)$ is
commutative up to homotopy \cite{G}.

If one wishes to furthermore relax the other conditions  ``up to
homotopy'', one can fix that $\eta$ needs to be non-degenerate
only on $HH^0(A,A)$ and  only require  that  $HH^0(A,A)$ has to be
finite dimensional. In this case, the operadic operations defined
below will give operations $f:A^{\otimes n}\rightarrow HH^0(A,A)$
and will thus give actions only up to homotopy. This is enough to
get the BV structure on $CH^*(A,A)$, but not quite enough to lift
the action to the chain level. We are currently working on such a
construction in formal geometry and defer the reader to
\cite{stringarc}.

\subsection{Notation} Let $(e_i)$ be a basis for $A$ and let $C:=
e_i \eta^{ij} \otimes e_j$ be the Casimir element, i.e.
$\eta^{ij}$ is the inverse to $\eta_{ij}=\eta(e_i,e_j)$.

With the help of the non--degenerate bilinear form, we identify
\begin{equation}
\label{identify} CH^n(A,A)= Hom(A^{\otimes n},A) \cong A\otimes
A^{* \otimes n} \cong A^{* \otimes n+1}
\end{equation}
 We would like to stress the order of the
tensor products we choose. This is the order from right to left,
which works in such a way that one does not need to permute tensor
factors in order to contract.

If $f \in Hom(A^{\otimes n},A)$, we denote by $\tilde f$ its image
in $A^{* \otimes n+1}$, explicitly $\tilde f(a_0, \dots
,a_n)=\eta(a_0,f(a_1,\dots,a_n))$.

With the help of (\ref{identify}) we can pull back the Connes'
operators $b$ and $B$ (see e.g.\ \cite{Loday}) on the spaces
$A^{\otimes n}$ to their duals and to $Hom(A^{\otimes n},A)$.

Also let $t:A^{\otimes n}\rightarrow A^{\otimes n}$ be the
operator given by performing a cyclic permutation $(a_1,\dots,a_n)
\mapsto (-1)^{n-1}(a_n,a_1,\dots a_{n-1})$ and $N:=1+t+ \cdots
+t^{n-1}:A^{\otimes n}\rightarrow A^{\otimes n}$.

 It is easy to check that the
 operator induced
by $b$ is exactly the Hochschild differential; we will denote this
operator by $\del$. We write $\Delta$ for the operator induced by
$B$. It follows that $\Delta^2=0$ and $\Delta\del+\del\Delta=0$.

\subsection{Assumption}
 To make the formulas simpler we will restrict to normalized
 Hochschild cochains $\overline {CH}^n(A,A)$ which are the $f\in
 CH^n(A,A)$ which vanish when evaluated on any tensor containing
 $1\in A$ as a tensor factor (see e.g.\ \cite{Loday}).
On the normalized chains the operator $\Delta$ is explicitly
defined as follows: for $f\in \overline {CH}^n(A,A)$
\begin{equation}
\eta(a_0,(\Delta f)(a_1,\dots a_{n-1})):=\eta(1,f\circ N(a_0,\dots
a_n))
\end{equation}

\subsection{Correlators from decorated trees}

We will use the notation of tensor products indexed by arbitrary
sets, see e.g.\ \cite{Deligne}. For a linearly ordered set $I$
 denote by $\bigcup_I a_i$ the product of the $a_i$ in the
order dictated by $I$.

\begin{df} Let $\t$ be the realization
of a spine decorated planted planar b/w tree, $v \in V_w$, and
$f\in \overline {CH}^{|v|}(A,A)$. We define
$Y(v,f):A^{F_v(\t)}\rightarrow k$ by

$$Y(v, f) (\bigotimes_{i\in F_v(\t)} a_i):=\eta(a_{\Fnum^{-1}(0)},
f(a_{\Fnum^{-1}(1)}\otimes \dots \otimes a_{\Fnum^{-1}(|v|)}))$$

Set $V_{b-int}:=V_b(\t)\setminus (V_{tail}\cup \{v_{root}\} \cup
V_{spine})$. For $v\in V_{b-int}$ we define $Y(v):=
A^{F_v(\t)}\rightarrow k$ by
 $$Y(v)(\bigoplus_{i\in
F_v(\t)} a_i)=\eta(1, \bigcup_{i\in F_{v}} a_i)$$

\end{df}

\begin{df}
\label{treeactionone} Let $\t$ be the realization of a planar
planted b/w tree with $n$ free tails and $k$ labels and $f_i \in
\overline{CH}^{n_i}(A,A)$. For such a tree there is a canonical
identification  $\{v_{root}\} \cup V_{ftail} \rightarrow
\{0,1,\dots,|V_{ftail}|\}$ which is given by sending $v_{root}$ to
$0$ and enumerating the tails in the linear order induced by the
planted planar tree. Set $E_{int}(\t):= E(\t)\setminus
(E_{tail}\cup E_{root} \cup E_{spine})$ and for
$(a_0,\dots,a_n)\in A^{\otimes (\{v_{root}\}\cup V_{ftail})}$  set
\begin{multline}
Y(\t)(f_1,\dots,f_k)(a_0,\dots,a_n) := \\
\left(\bigotimes_{v\in V_w(\t)} Y(v,f_{\lab(v)})\bigotimes_{v\in
V_{b-int}} Y_v\right)\left( (\bigotimes_{i\in F_{ftail}(\t)\cup
\{F_{root}\}}a_i)(\bigotimes_{j\in F_{spine}} 1) \otimes
C^{\otimes E_{int}(\t)}\right)\end{multline}
\end{df}

In other words,  decorate the root flag by $a_0$, the free tail
flags by $a_1,\dots,a_n$, the spines by $1$ and the edges by $C$
and then contract tensors according to the decoration at the white
vertices while using the product at the black vertices.

\begin{df}
\label{degreecount} We extend the definition above by
\begin{equation}
Y(\t)(f_1,\dots,f_k)(a_0,\dots,a_n)=0 \text{ if }
|v_{\lab^{-1}(i)}| \neq n_i=:|f_i|
\end{equation}

\end{df}

\subsection{The foliage operator} Let $F$ be the foliage operator of \cite{del}
applied to trees. This means that $F(\t)$ is the formal sum over
all trees obtained from $\t$ by gluing an arbitrary number of free
tails to the white vertices. The extra edges are called free tail
edges $E_{ftail}$ and the extra vertices $V_{ftail}$ are defined
to be black and are called free tail vertices.

Using the trees defined in Figure \ref{cactexamples} this
corresponds to the formal sum $F(\t):= \sum_n l_n \circ_v \t$
where the operadic composition is the one for b/w trees which are
not necessarily bi-partite (see \cite{del}). In our current setup
we should first form $\tilde F(\t):=\sum_n \t_n \circ_v \t$ and
then delete the images of all leaf edges together with their white
vertices of the $\t_n$ to obtain $F(\t)$.

\subsection{Signs}
\label{signsection}  The best way to fix signs of course is to
work with tensors indexed by edges like in \cite{del,KS}. For this
one fixes a free object $L$ (free $\mathbb{Z}$-module or
$k$-vector space) generated by one element of degree $\pm 1$ and
calculates signs using $L^{\otimes E_w(\t)}$ before applying the
foliage operator while using $L^{\otimes E_{weight}}$  after
applying the foliage operator, where $E_{weight}=E_w\cup
E_{root}\cup E_{ftail}\cup E_{spine}$.

Explicitly, we fix the signs to be given as follows. For any tree
$\t'$ in the linear combination above, we take the sign of $\t'$
to be the sign of the permutation which permutes the set
$E_{weight}$ in the order induced by $\prec$  to the order where
at each vertex one first has the root if applicable, then all
non--tail edges, then all the free tails, and if there is a spine
edge, the spine.

The explicit signs above coincide with usual signs  \cite{Loday}
for the operations and the operators $b$ and $B$ and also coincide
with the signs of \cite{G} for the $\circ_i$ and hence for the
brace operations. The signs for the operations corresponding to
operations on the Hochschild side are fixed  by declaring the
symbols ``,'' and ``\{'' to have degree one.

\begin{df}
\label{treeaction}
 For $\t\in \swlbptree$ let $\hat\t$ be its realization. We define the
operation of $\t$ on $\overline {CH}(A,A)$ by
\begin{equation}
\eta(a_0,{\t(f_1,\dots, f_n)}(a_1,\dots,
a_N)):=Y(F(\hat\t))(f_1,\dots, f_n)(a_0,\dots,a_N)
\end{equation}
Notice that due to the Definition \ref{degreecount} the right hand
side is finite.
\end{df}

\subsection{Examples}
\label{example}
 We will first regard the tree $O'$ with one white vertex, no
additional black edges, no free tails and a spine, see Figure
\ref{cactexamples}. For a function $f\in \overline {CH}^n$ we
obtain:
\begin{multline*}
Y(F(O'))(f)(a_0,\dots,a_{n-1})= \eta(1,f(a_0,\dots
a_{n-1})+(-1)^{n-1}f(a_{n-1}, a_0, \dots, a_{n-2}) +\dots)\\
=\eta(a_0,\Delta(f)(a_1,\dots,a_{n-1})) \end{multline*}

Let $\t'_{n,i}$ be the tree of Figure \ref{cactexamples}. Then the
operation corresponds to
$$
Y(F(\t'_{n,i}))(f;g_1,\dots,g_n)(a_0,\dots,a_{N})=\\
\eta(1,f\{' g_{i+1}, \dots, g_n, g_1, \dots, g_i\}(a_{(2)},
a_0,a_{(1)}))
$$
where $N=|f|+\sum|g_i|-n-1$ and we used the short hand notation
\begin{multline*}
f\{' g_{j+1}, \dots, g_n, g_1, \dots, g_j\}(a_{(2)}, a_0,a_{(1)}) =
\sum \pm
f(a_{k+1},\dots, a_{i_{j+1}-1},\\ g_{j+1}(a_{i_{j+1}}, \dots,
a_{i_{j+1}+|g_{j+1}|}), \dots, a_{i_n-1}, g_n(a_{i_n}, \dots,
   a_{i_n+|g_n|}), \dots , a_N,
   a_0, \\a_1, \dots, a_{i_1-1},g_1(a_{i_1}, \dots, a_{i_1+|g_1|}), \dots,
a_{i_j-1}, g_j(a_{i_j}, \dots, a_{i_j+|g_j|}), \dots , a_k)
\end{multline*}
where the sum runs over
$1 \leq i_1 \leq  \dots \leq i_j
\leq  \dots \leq k \leq  \dots \leq i_{j+1} \leq  \dots \leq i_{n} \leq N:$
 $i_l + |g_l|\leq i_{l+1},i_j + |g_j|\leq k $ and the signs are
as explained above.

\begin{thm}[The cyclic Deligne conjecture]
\label{mainthm}
The Hochschild cochains of a finite-dimensional associative algebra with a
non--degenerate, symmetric, invariant, bilinear form are an
algebra over the chains of the framed little discs operad. This
operation is compatible with the differentials.
\end{thm}

\begin{proof}
We will use the cellular chains $\CCCacti$ as a model for the chains
of the framed little discs operad. It is clear that
\ref{treeaction} defines an action.
 On the Hochschild side, the $\circ_i$ operations
are substitutions of the type $f_i=\psi(g_1,\dots,g_n)$. For
$\CCCacti$ the $\t \circ_i \t'$ operations are the pull-back via
the foliage operator of all possible substitutions of elements of
$F(\t), \t \in \CCCacti$ into the position $i$ of $F(\t'$). The
action  $Y$ then projects onto the substitution
$f_i=\psi(g_1,\dots,g_n)$ so that the action is operadic.
Explicitly the substitution $t \circ^s_i t'$ for planted planar
bi-partite trees with a decoration $\sdec$ and additional free
tails is given as follows: Say the number of tails of $t'$
coincides with $|F(v_i)|$. In this case replace the vertex $v_i$
of $t$, its edges and the black vertices corresponding to the
edges with the tree $t'$ matching the flags of $v_i$ with the
tails of $t'$  by first matching the root edge with the marked
flag of $v_i$ and then using the linear order. Lastly contract the
image of the root flag. Otherwise set $t \circ^s_i t'=0$. With
this definition it is easy to see that $F(\t \circ \t')=F(\t)
\circ^s_i F(\t')$.

The compatibility of the Hochschild differential with the
differential of the cell complex follows from the relevant
statements for $\t_n$ and $\t_n^b$, which are a straightforward
but lengthy calculation (see e.g.\ \cite{del,G}), together with
the calculations above \S\ref{example} which are easily modified
to show that $(\del O')(f)=\Delta(\del(f))$ and that
$(\del\t'_{n,i})(f,g_1, \dots, g_n)= (\del\t'_{n,i})(f,g_1, \dots,
g_n) \pm (\t'_{n,i})(\del f,g_1, \dots, g_n) + \sum_i\pm
 (\t'_{n,i})(f,g_1, \dots, \del(g_i), \dots, g_n)$
via an even more lengthy but still straightforward calculation.
This then verifies the claim in view of the compatibility of the
differentials and the respective operad structures.

Alternatively,  in view of the operation of the foliage operator,
the compatibilities follow from a straightforward translation of
trees with tails into operations on the Hochschild complex. The
compatibility of the differential then follows from the almost
identical  definition of the differential for trees with tails of
\cite{del} and that in the Hochschild complex as $\del(f)=f\circ
\cup - (-1)^{|f||}\cup \circ f$.
\end{proof}

\begin{cor}
The normalized Hochschild cochains of an algebra as above are a BV
algebra up to homotopy.
\end{cor}

This could of course have been checked directly without recourse to the
operation of a chain model, but we do not know of any source for
this result. It also seems to be difficult to guess the right
 homotopies as Gerstenhaber did in the non-cyclic case \cite{G}.
The content of the next corollary was expected
\cite{Connes}, but we again
could not find a source for it.

\begin{cor}
The Hochschild cohomology of an algebra as above is a BV algebra,
such that the induced bracket is the Gerstenhaber bracket.
\end{cor}

Lastly, since our second version of cellular chains of Proposition
\ref{secondcells} are a subdivision of the cell decomposition of
Proposition \ref{firstcells}, we can also use the latter cell
decomposition.

\begin{cor}
The normalized Hochschild cochains of an algebra as above are an
algebra over the semi--direct product over a chain model of the
little discs operad and a chain model for the operad $\mathcal{S}$
built on the monoid $S^1$.
\end{cor}

\begin{rmk}
The operation of the little discs operad by braces, viz.\ the
original Deligne conjecture as discussed in \cite{del} for
Frobenius algebras, corresponds to the decorations in which $\Zdec
\equiv 0$ and the decorated edge is always the outgoing edge.
\end{rmk}

\begin{rmk}
In the Theorem \ref{mainthm} we can relax the conditions and
implications as explained in \S\ref{assumption}.
\end{rmk}

\section{Variations and relation to string topology}
In terms of the setup of operadic correlation functions which we
presented above, it is possible to analyze several
generalizations. First, one can generalize from trees to more
general graphs. This description then yields an action of the
pseudo-cells of moduli spaces of curves or bordered surfaces
\cite{ribbon}. One can also consider different types of chains,
such as Hochschild chains or cyclic (co)--chains. The latter also
works well with omitting markings to the trees or regarding
unmarked graphs \cite{ribbon}.

In \cite{KLP} we gave a map called loop which maps the so--called
$\Arc$ operad to ribbon graphs with marked points on the cycles of
the graph. In the case of no punctures the analysis of this map in
terms of Strebel differentials yields another proof of Penner's
theorem \cite{P} on the homotopy equivalence of the suboperad of
quasi--filling arcs and the moduli space of decorated bordered
surfaces  \cite{ribbon}. This in turn gives a cell decomposition
of the aforementioned moduli space. Moreover the correspondence
induces an operadic structure on ribbon graphs by pulling back the
gluings from the $\Arc$ operad. Using the operadic correlation
functions it is straightforward to obtain an action of the cells
on a cyclic complex. In a similar spirit, an action of the framed
little discs on a cyclic complex given by the $Tot$ of a special
type of cyclic cosimplicial complex has been announced in
\cite{MS}. Constructing the action  in terms of our correlation
functions then should allow us to construct an operation of the
cells of moduli space on such a complex. Moreover a further
decoration of the cells by $\mathbb{Z}/2\mathbb{Z}$ produces an
operad which acts on the cyclic complex of such an algebra and is
compatible with the differential \cite{ribbon}.

The $A_{\infty}$-- versions of these statements could be deduced
from a conjectural ``blow--up''  of the cacti operads which is
presented in
 \cite{woods}. Here the cells are given by products of
 associahedra and cyclohedra
and are indexed by trees of the type appearing in \cite{KS}.

Finally using the cyclic description of the free loop space or the
iterated integral representation of \cite{Merk} together with the
results mentioned above, we expect to be able to obtain an action
of the (decorated) pseudo-cells of moduli space action on the free
loop space of a compact manifold which extends the operation of
the string PROP \cite{CS1,CS2} thus completing a further step of
the string topology program \cite{stringarc}.

\end{document}